\documentclass[12pt]{amsart}
\usepackage{amssymb}
\usepackage[all]{xy}

\theoremstyle{plain}
\newtheorem{theorem}{Theorem}[section]
\newtheorem{corollary}[theorem]{Corollary}
\newtheorem{lemma}[theorem]{Lemma}
\newtheorem{proposition}[theorem]{Proposition}
\newtheorem{conjecture}[theorem]{Conjecture}

\theoremstyle{definition}

\newtheorem{remark}[theorem]{Remark}
\newtheorem*{question}{Question}

\theoremstyle{remark}

\newcommand{\abs}[1]{\lvert#1\rvert}
\newcommand{\norm}[1]{\lVert#1\rVert}

\newcommand{\bignorm}[1]{\bigl\lVert#1\bigr\rVert}

\renewcommand{\le}{\leqslant}
\renewcommand{\ge}{\geqslant}
\renewcommand{\mid}{\::\:}

\newcommand{\term}[1]{{\textit{\textbf{#1}}}}

\DeclareMathOperator{\supp}{supp}
\DeclareMathOperator{\Range}{Range}
\DeclareMathOperator{\Span}{span}
\DeclareMathOperator{\sign}{sign}

\begin{document}
\baselineskip 18pt

\title[Operator ideals in Lorentz spaces]
      {Norm closed operator ideals\\ in Lorentz sequence spaces}

\author[A.~Kaminska]{Anna Kaminska}
\author[A.~I.~Popov]{Alexey I. Popov}
\author[E.~Spinu]{Eugeniu Spinu}
\author[A.~Tcaciuc]{Adi Tcaciuc}
\author[V.~G.~Troitsky]{Vladimir G. Troitsky}

\address[A.~Kaminska]
  {Department of Mathematical Sciences
   The University of Memphis, Memphis, TN 38152-3240. USA}
\email{kaminska@memphis.edu}

\address[A.~I.~Popov, E.~Spinu, and V.~G.~Troitsky]
  {Department of Mathematical
  and Statistical Sciences, University of Alberta, Edmonton,
  AB, T6G\,2G1. Canada}
\email{apopov@math.ualberta.ca, spinu@ualberta.ca}
\email{troitsky@ualberta.ca}

\address[A.~Tcaciuc]{Mathematics and Statistics Department,
   Grant MacEwan University, Edmonton, AB, T5J\,P2P, Canada}
\email{atcaciuc@ualberta.ca}

\thanks{The forth and fifth authors were supported by NSERC}

\date{\today}

\begin{abstract}
In this paper, we study the structure of closed algebraic ideals in the algebra of operators acting on a Lorentz sequence space. 
\end{abstract}

\maketitle

\section{Introduction}\label{introduction}

\subsection{Ideals} This paper is concerned with the study of the structure of closed algebraic ideals in the algebra $L(X)$ of all bounded linear operatrors on a Banach space~$X$. 

Throughout the paper, by a \term{subspace} of a Banach space we mean a closed subspace; a vector subspace of $X$ which is not necessarily closed will be referred to as \term{linear subspace}. A (two-sided) \term{ideal} in $L(X)$ is a linear subspace $J$ of $L(X)$ such that $ATB\in J$ whenever $T\in J$ and $A,B\in L(X)$. The ideal $J$ is called \term{proper} if $J\ne L(X)$. The ideal $J$ is \term{non-trivial} if $J$ is proper and $J\ne\{0\}$.

The spaces for which the structure of closed ideals in $L(X)$ is well-understood are very few. It was shown in~\cite{Calk} that the only non-trivial closed ideal in the algebra $L(\ell_2)$ is the ideal of compact operators. This result was generalized in~\cite{FGM} to the spaces $\ell_p$ ($1\le p<\infty$) and~$c_0$. A space constructed recently in~\cite{AH} is another space with this property. In~\cite{LLR04} and~\cite{LSZ06}, it was shown that the algebras $L\big((\oplus_{k=1}^\infty\ell_2^k)_{c_0}\big)$ and $L\big((\oplus_{k=1}^\infty\ell_2^k)_{\ell_1}\big)$ have exactly two non-trivial closed ideals. There are no other separable spaces for which the structure of closed ideals in $L(X)$ is completely known.

Partial results about the structure of closed ideals in $L(X)$ were obtained in~\cite[5.3.9]{Pie81} for $X=L_p[0,1]$ ($1<p<\infty$, $p\ne 2$) and in~\cite{SSTT} and~\cite{Schl11} for $L(\ell_p\oplus\ell_q)$ ($1\le p,q<\infty$). The purpose of this paper is to investigate the structure of ideals in $L(d_{w,p})$ where $d_{w,p}$ is a Lorentz sequence space (see the definition in Subsection~\ref{lorentz-subsec}).

For two closed ideals $J_1$ and $J_2$ in $L(X)$, we will denote by $J_1\wedge J_2$ the largest closed ideal $J$ in $L(X)$ such that $J\subseteq J_1$ and $J\subseteq J_2$ (that is, $J_1\wedge J_2=J_1\cap J_2$), and we will denote by $J_1\vee J_2$ the smallest closed ideal $J$ in $L(X)$ such that $J_1\subseteq J$ and $J_2\subseteq J$. We say that $J_2$ is a \term{successor} of $J_1$ if $J_1\subsetneq J_2$. If, in addition, no closed ideal $J$ in $L(X)$ satisfies $J_1\subsetneq J\subsetneq J_2$, then we call $J_2$ an \term{immediate successor} of~$J_1$. 

It is well-known that if $X$ is a Banach space then every non-zero ideal in the algebra $L(X)$ must contain the ideal $\mathcal F(X)$ of all finite-rank operators on~$X$. It follows that, at least in the presence of the approximation property (in particular, if $X$ has a Schauder basis), every non-zero closed ideal in $L(X)$ contains the closed ideal $\mathcal K(X)$ of all compact operators. 

Two ideals closely related to $\mathcal K(X)$ are the closed ideal $\mathcal{SS}(X)$ of strictly singular operators and the closed ideal $\mathcal{FSS}(X)$ of finitely strictly singular operators on~$X$. Recall that an operator $T\in L(X)$ is called \term{strictly singular} if no restriction $T|_{Z}$ of $T$ to an infinite-dimensional subspace $Z$ of $X$ is an isomorphism. An operator $T$ is \term{finitely strictly singular} if for any $\varepsilon>0$ there is $N\in\mathbb N$ such that any subspace $Z$ of $X$ with $\dim Z\ge N$ contains a vector $z\in Z$ satisfying $\norm{Tz}<\varepsilon\norm{z}$. It is not hard to show that $\mathcal K(X)\subseteq\mathcal{FSS}(X)\subseteq{SS}(X)$ (see \cite{LT77,Mil70,SSTT,ADST} for more information about these classes of operators).

If $X$ is a Banach space and $T\in L(X)$ then the ideal in $L(X)$ generated by $T$ is denoted by~$J_T$. It is easy to see that $J_T=\bigl\{\sum_{i=1}^nA_iTB_i\mid A_i,B_i\in L(X)\bigr\}$. It follows that if $S\in L(X)$ factors through~$T$, i.e., $S=ATB$ for some $A,B\in L(X)$ then $J_S\subseteq J_T$.

\subsection{Basic sequences} The main tool in this paper is the notion of a basic sequence. In this subsection, we will fix some terminology and remind some classical facts about basic sequences. For a thorough introduction to this topic, we refer the reader to~\cite{Car} or~\cite{Zizler}. 

If $(x_n)$ is a sequence in a Banach space $X$ then its closed span will be denoted by $[x_n]$. We say that a basic sequence $(x_n)$ \term{dominates} a basic sequence $(y_n)$ and write $(x_n)\succeq(y_n)$ if the convergence of a series $\sum_{n=1}^\infty a_nx_n$ implies the convergence of the series $\sum_{n=1}^\infty a_ny_n$. We say that $(x_n)$ is \term{equivalent} to $(y_n)$ and write $(x_n)\sim(y_n)$ if $(x_n)\succeq(y_n)$ and $(y_n)\succeq(x_n)$.

\begin{remark}
It follows from the Closed Graph Theorem that $(x_n)\succeq(y_n)$ if and only if the linear map from $\Span\{x_n\}$ to $\Span\{y_n\}$ defined by the formula $T\colon x_n\mapsto y_n$ is bounded.
\end{remark}

If $(x_n)$ is a basis in a Banach space~$X$, $z=\sum_{i=1}^\infty z_ix_i\in X$, and $A\subseteq\mathbb N$ then the vector $\sum_{i\in A} z_ix_i$ will be denoted by $z|_A$ (provided the series converges; this is always the case when the basis is unconditional). We will refer to $z|_A$ as the \term{restriction of $z$ to $A$}. The restrictions $z|_{[n,\infty)\cap\mathbb N}$ and $z|_{(n,\infty)\cap\mathbb N}$, where $n\in\mathbb N$, will be abbreviated as $z|_{[n,\infty)}$ and $z|_{(n,\infty)}$, respectively. We say that a vector $v$ is a \term{restriction} of $z$ if there exists $A\subseteq\mathbb N$ such that $v=z|_A$. The vector $z=\sum_{i=1}^\infty z_ix_i$ will also be denoted by $z=(z_i)$. If $z=\sum_{i=1}^\infty z_ix_i$ then the \term{support} of $z$ is the set $\supp z=\{i\in\mathbb N\colon z_i\ne 0\}$.

Every 1-unconditional basis $(x_n)$ in a Banach space $X$ defines a Banach lattice order on $X$ by $\sum_{i=1}^\infty a_ix_i\ge 0$ if and only if $a_i\ge 0$ for all $i\in\mathbb N$ (see, e.g., \cite[page~2]{LT79}). For $x\in X$, we have $\abs{x}=x\vee(-x)$. A Banach lattice is said to have \term{order continuous norm} if the condition $x_\alpha\downarrow 0$ implies $\norm{x_\alpha}\to 0$. For an introduction to Banach lattices and standard terminology, we refer the reader to~\cite[\S1.2]{AA02}.

If $(x_n)$ is a basic sequence in a Banach space~$X$, then a sequence $(y_n)$ in $\Span\{x_n\}$ is a \term{block sequence} of $(x_n)$ if there is a strictly increasing sequence $(p_n)$ in $\mathbb N$ and a sequence of scalars $(a_i)$ such that $y_n=\sum_{i=p_{n}+1}^{p_{n+1}}a_ix_i$ for all $n\in\mathbb N$.

The following two facts are classical and will sometimes be used without any references. The first fact is known as the Principle of Small Perturbations (see, e.g., \cite[Theorem~4.23]{Zizler}).

\begin{theorem}\label{small-perturbations}
Let $X$ be a Banach space, $(x_n)$ a basic sequence in~$X$, and $(x_n^*)$ the correspondent biorthogonal functionals defined on $[x_n]$. If $(y_n)$ is a sequence such that $\sum_{n=1}^\infty\norm{x^*_n}\cdot\norm{x_n-y_n}<1$ then $(y_n)$ is a basic sequence equivalent to $(x_n)$. Moreover, if $[x_n]$ is complemented in $X$ then so is $[y_n]$. If $[x_n]=X$ then $[y_n]=X$.
\end{theorem}

The next fact, which is often called the Bessaga-Pe\l czy\'nski selection principle, is a result of combining the ``gliding hump'' argument (see, e.g., \cite[Lemma~5.1]{Car}) with the Principle of Small Perturbations. 

\begin{theorem}\label{bessaga-pelczynski}
Let $X$ be a Banach space with a seminormalized basis $(x_n)$ and let $(x_n^*)$ be the correspondent biorthogonal functionals. Let $(y_n)$ be a seminormalized sequence in $X$ such that $x^*_n(y_k)\stackrel{k\to\infty}{\longrightarrow} 0$ for all $n\in\mathbb N$. Then $(y_n)$ has a subsequence $(y_{n_k})$ which is basic and equivalent to a block sequence $(u_k)$ of $(x_n)$. Moreover, $y_{n_k}-u_k\to 0$, and $u_k$ is a restriction of $y_{n_k}$.
\end{theorem}

\subsection{Lorentz sequence spaces}\label{lorentz-subsec}

Let $1\le p<\infty$ and $w=(w_n)$ be a sequence in $\mathbb R$ such that $w_1=1$, $w_n\downarrow 0$, and $\sum_{i=1}^\infty w_i=\infty$. The Lorentz sequence space $d_{w,p}$ is a Banach space of all vectors $x\in c_0$ such that $\norm{x}_{d_{w,p}}<\infty$, where 
$$
\norm{(x_n)}_{d_{w,p}}=\Big(\sum_{n=1}^\infty w_nx_n^{*p}\Big)^{1/p}
$$
is the norm in~$d_{w,p}$. Here $(x_n^*)$ is the \term{non-increasing rearrangement} of the sequence $(\abs{x_n})$. An overview of  properties of Lorentz sequence spaces can be found in \cite[Section~4.e]{LT77}.

The vectors $(e_n)$ in $d_{w,p}$ defined by $e_n(i)=\delta_{ni}$ ($n,i\in\mathbb N$) form a 1-symmetric basis in~$d_{w,p}$. In particular, $(e_n)$ is 1-unconditional, hence $d_{w,p}$ is a Banach lattice. We call $(e_n)$ the unit vector basis of~$d_{w,p}$. The unit vector basis of $\ell_p$ will be denoted by $(f_n)$ throughout the paper. 

\begin{remark}\label{ell-p-wealth}
It is proved in \cite[Lemma~1]{ACL} and \cite[Lemma~15]{CL} that if $(u_n)$ is a seminormalized block sequence of $(e_n)$ in~$d_{w,p}$, $u_n=\sum_{i=p_n+1}^{p_{n+1}}a_ie_i$, such that $a_i\to 0$, then there is a subsequence $(u_{n_k})$ such that $(u_{n_k})\sim(f_n)$ and $[u_{n_k}]$ is complemented in~$d_{w,p}$. Further, it was shown in \cite[Corollary~3]{ACL} that if $(y_n)$ is a seminormalized block sequence of $(e_n)$ then there is a seminormalized block sequence $(u_n)$ of $(y_n)$ such that $u_n=\sum_{i=p_n+1}^{p_{n+1}}a_ie_i$, with $a_i\to 0$. Therefore, every infinite dimensional subspace of $d_{w,p}$ contains a further subspace which is complemented in $d_{w,p}$ and isomorphic to $\ell_p$ (\cite[Corollary~17]{CL}).
\end{remark}

\begin{remark}\label{wsc-ocn}
Remark~\ref{ell-p-wealth} yields, in particular, that $d_{w,p}$ does not contain copies of~$c_0$. Since the basis $(e_n)$ of $d_{w,p}$ is unconditional, the space $d_{w,p}$ is weakly sequentially complete by \cite[Theorem~4.60]{AB} (see also \cite[Theorem 1.c.10]{LT77}). Also, \cite[Theorem~4.56]{AB} guarantees that $d_{w,p}$ has order continuous norm. In particular, if $x\in d_{w,p}$ then $\norm{x|_{[n,\infty)}}\to 0$ as $n\to\infty$.
\end{remark}

\begin{remark}
It was shown in \cite{G69} that if $p>1$ then $d_{w,p}$ is reflexive. This can also be easily obtained from Remark~\ref{ell-p-wealth} (cf. \cite[Theorem 1.c.12]{LT77}).
\end{remark}

\begin{remark}\label{e-n-w-null-remark}
The unit vector basis $(e_n)$ of $d_{w,p}$ is weakly null. Indeed, by Rosenthal's $\ell_1$-theorem (see~\cite{Ros74}; also \cite[Theorem~2.e.5]{LT77}), $(e_n)$ is weakly Cauchy. Since it is symmetric, $(e_n)\sim(e_{2n}-e_{2n-1})$.
\end{remark}

The next proposition will be used often in this paper.

\begin{proposition}[{\rm\cite[Proposition~5 and Corollary~2]{ACL}}]\label{domination}
If $(u_n)$ is a seminormalized block sequence of $(e_n)$ then $(f_n)\succeq(u_n)$. If $(u_n)$ does not contain subsequences equivalent to $(f_n)$ then also $(u_n)\succeq(e_n)$.
\end{proposition}

The following lemma is standard.

\begin{lemma}\label{lattice-domination}
Let $(x_n)$ be a block sequence of $(e_n)$, $x_n=\sum_{i=p_n+1}^{p_{n+1}}a_ie_i$. If $(y_n)$ is a basic sequence such that $y_n=\sum_{i=p_n+1}^{p_{n+1}}b_ie_i$, where $\abs{b_i}\le\abs{a_i}$ for all $i\in\mathbb N$, then $(x_n)$ is basic and $(x_n)\succeq(y_n)$.
\end{lemma}
\begin{proof}
Let 
$$
\gamma_i=\left\{\begin{array}{ll}
\frac{b_i}{a_i},\quad&\mbox{if }a_i\ne 0,\\
0,&\mbox{if }a_i=0.
\end{array}\right.
$$
Define an operator $T\in L(d_{w,p})$ by $T\left(\sum_{i=1}^\infty c_ie_i\right)=\sum_{i=1}^\infty c_i\gamma_ie_i$. Then $T$ is, clearly, linear and, since the basis $(e_n)$ is 1-unconditional, $T$ is bounded with $\norm{T}\le 1$. In particular, $T|_{[x_n]}$ is bounded. Also, $T(x_n)=y_n$ for all $n\in\mathbb N$, hence $(x_n)\succeq(y_n)$.
\end{proof}

\subsection{Outline of the results}

The purpose of the paper is to uncover the structure of ideals in $L(d_{w,p})$. We show that (some of) these ideals can be arranged into the following diagram.
$$
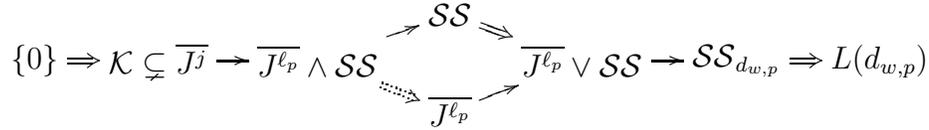
\xymatrix @R=0pt@C=12pt{ & &  &\mathcal{SS}\ar@{=>}[dr]& & & \\\{0\}\ar@{=>}[r]& \mathcal K\subsetneq\overline{J^j}\ar[r]&\overline{J^{\ell_p}}\wedge\mathcal{SS}\ar@{:>}[dr]\ar[ur]& &\overline{J^{\ell_p}}\vee\mathcal{SS}\ar[r]&\mathcal{SS}_{d_{w,p}}\ar@{=>}[r]&L(d_{w,p})\\ & & &\overline{J^{\ell_p}}\ar[ur]& & & }
$$
(the notations will be defined throughout the paper). On this diagram, a single arrow between ideals, $J_1\longrightarrow J_2$, means that $J_1\subseteq J_2$. A double arrow between ideals, $J_1\Longrightarrow J_2$, means that $J_2$ is the only immediate successor of $J_1$ (in particular, $J_1\ne J_2$), whereas a dotted double arrow between ideals, $\xymatrix @C=20pt{J_1\ar@{:>}[r]&J_2}$, only shows that $J_2$ is an immediate successor for $J_1$ (in particular, $J_1$ may have other immediate successors). 

While working with the diagram above, we obtain several important characterizations of some ideals in $L(d_{w,p})$. In particular, we show that $\mathcal{FSS}(d_{w,p})=\mathcal{SS}(d_{w,p})$ (Theorem~\ref{FSS}). We also characterize the ideal of weakly compact operators (Theorem~\ref{weakly-compact}) and Dunford-Pettis operators (Theorem~\ref{dunford-pettis}) on~$d_{w,p}$. We show in Theorem~\ref{J-j-direct} that $\overline{J^j}$ is the only immediate successor of $\mathcal K$ under some assumption on the weights~$w$. In the last section of the paper, we show that all strictly singular operators from $\ell_1$ to $d_{w,1}$ can be approximated by operators factoring through the formal identity operator $j\colon\ell_1\to d_{w,1}$ (see Section~\ref{form-dentity} for the definition). We also obtain a result on factoring positive operators from $\mathcal{SS}(d_{w,p})$ through the formal identity operator (Theorem~\ref{positive-factor}).

\section{Operators factorable through $\ell_p$}\label{ell-p-section}

Let $X$ and $Y$ be Banach spaces and $T\in L(X)$. We say that \term{$T$ factors through $Y$} if there are two operators $A\in L(X,Y)$ and $B\in L(Y,X)$ such that $T=BA$. 

The following two lemmas are standard. We present their proofs for the sake of completeness.

\begin{lemma}\label{ST=I}
Let $X$ and $Y$ be Banach spaces and $T\in L(X,Y)$, $S\in L(Y,X)$ be such that $ST=\mathrm{id}_X$. Then $T$ is an isomorphism and $\Range T$ is a complemented subspace of $Y$ isomorphic to~$X$.
\end{lemma}

\begin{proof}
For all $x\in X$, we have $\norm{x}=\norm{STx}\le\norm{S}\norm{Tx}$, so $\norm{Tx}\ge\frac{1}{\norm{S}}\norm{x}$. This shows that $T$ is an isomorphism. In particular, $\Range T$ is a closed subspace of $Y$ isomorphic to~$X$.

Put $P=TS\in L(Y)$. Then $P^2=TSTS=T\mathrm{id}_XS=TS=P$, hence $P$ is a projection. Clearly, $\Range P\subseteq\Range T$. Also, $PT=TST=T$, so $\Range T\subseteq\Range P$. Therefore $\Range P=\Range T$, and $\Range T$ is complemented.
\end{proof}

\begin{lemma}\label{gen-fact-ideal}
Let $X$ and $Y$ be Banach spaces such that $Y$ is isomorphic to $Y\oplus Y$. Then the set $J=\bigl\{T\in L(X)\mid T\mbox{ factors through }Y\bigr\}$ is an ideal in $L(X)$. 
\end{lemma}

\begin{proof}
It is clear that $J$ is closed under multiplication by operators in  $L(X)$. In particular, $J$ is closed under scalar multiplication. Let $A,B\in J$.  Write $A=A_1A_2$ and $B=B_1B_2$, where $A_1,B_1\in L(Y,X)$ and $A_2,B_2\in L(X,Y)$. Then $A+B=UV$ where $V\colon x\in X\mapsto(A_2x,B_2x)\in Y\oplus Y$ and $U\colon(x,y)\in Y\oplus Y\mapsto A_1x+B_1y\in Y$. Clearly, $UV$ factors through $Y\oplus Y\simeq Y$. Hence $A+B\in J$.
\end{proof}

We will denote the set of all operators in $L(d_{w,p})$ which factor through a Banach space $Y$ by~$J^Y$.

\begin{theorem}\label{ell-p-factor}
The sets $J^{\ell_p}$ and $\overline{J^{\ell_p}}$ are proper ideals in $L(d_{w,p})$.
\end{theorem}

\begin{proof}
Since $\ell_p$ is isomorphic to $\ell_p\oplus\ell_p$, it follows from Lemma~\ref{gen-fact-ideal} that $J^{\ell_p}$ is an ideal in $L(d_{w,p})$. Let us show that $J^{\ell_p}\ne L(d_{w,p})$.

Assume that $J^{\ell_p}=L(d_{w,p})$, then the identity operator $I$ on $d_{w,p}$ belongs to~$J$. Write $I=ST$ where $T\in L(d_{w,p},\ell_p)$ and $S\in L(\ell_p,d_{w,p})$. By Lemma~\ref{ST=I}, the range of $T$ is complemented in $\ell_p$ and is isomorphic to~$d_{w,p}$. This is a contradiction because all complemented infinite-dimensional subspaces of $\ell_p$ are isomorphic to $\ell_p$ (see, e.g., \cite[Theorem 2.a.3]{LT77}), while $d_{w,p}$ is not isomorphic to $\ell_p$ (see~\cite{CH} for the case $p=1$ and~\cite{G69} for the case $1<p<\infty$; see also \cite[p.~176]{LT77}).

Being the closure of a proper ideal, $\overline{J^{\ell_p}}$ is itself a proper ideal (see, e.g., \cite[Corollary VII.2.4]{Conway}).
\end{proof}

\begin{proposition}\label{ell-p-fact-characterization}
There exists a projection $P\in L(d_{w,p})$ such that $\Range P$ is isomorphic to~$\ell_p$. For every such $P$ we have $J_P=J^{\ell_p}$.
\end{proposition}

\begin{proof}
Such projections exist by Remark~\ref{ell-p-wealth}. Let $Y=\Range P$, $U\colon Y\to\ell_p$ be an isomorphism onto, and $i\colon Y\to d_{w,p}$ be the inclusion map. It is easy to see that $P=(iU^{-1})(UP)$, hence $P\in J^{\ell_p}$, so that $J_P\subseteq J^{\ell_p}$.

On the other hand, if $T\in J^{\ell_p}$ is arbitrary, $T=AB$ with $A\in L(\ell_p,d_{w,p})$, $B\in L(d_{w,p},\ell_p)$, then one can write $T=(AUP)P(iU^{-1}B)$, so that $T\in J_P$. Thus $J^{\ell_p}\subseteq J_P$.
\end{proof}

\begin{corollary}\label{K-ne-J-ell-p}
The ideal $\overline{J^{\ell_p}}$ properly contains the ideal of compact operators $\mathcal K(d_{w,p})$.
\end{corollary}

\begin{proof}
It was already mentioned in the introductory section that compact operators form the smallest closed ideal in $L(d_{w,p})$. Since a projection onto a subspace isomorphic to $\ell_p$ is not compact, it follows that $\mathcal K(d_{w,p})\ne\overline{J^{\ell_p}}$.
\end{proof}

\section{Strictly singular operators}

In this section we will study properties of strictly singular operators in $L(d_{w,p})$. Since projections onto the subspaces of $d_{w,p}$ isomorphic to $\ell_p$ are clearly not strictly singular, it follows from Proposition~\ref{ell-p-fact-characterization} that $\mathcal{SS}(d_{w,p})\ne J^{\ell_p}$. Moreover, $\mathcal{SS}\ne\overline{J^{\ell_p}}\vee\mathcal{SS}$ and $\overline{J^{\ell_p}}\wedge\mathcal{SS}\ne\overline{J^{\ell_p}}$. So, the ideals we discussed so far can be arranged as follows:
$$
\xymatrix @R=0pt@C=20pt{ & & &\mathcal{SS}\ar^{\ne}[dr]& & \\\{0\}\ar@{=>}[r]& \mathcal K\ar[r]&\overline{J^{\ell_p}}\wedge\mathcal{SS}\ar_{\ne}[dr]\ar[ur]& &\overline{J^{\ell_p}}\vee\mathcal{SS}\ar[r]&L(d_{w,p})\\ & & &\overline{J^{\ell_p}}\ar[ur]& & }
$$

The following theorem shows that there can be no other closed ideals between $\mathcal{SS}$ and $\overline{J^{\ell_p}}\vee\mathcal{SS}$ on this diagram.

\begin{theorem}\label{J-P-V-SS}
Let $T\in L(d_{w,p})$. If $T\not\in\mathcal{SS}(d_{w,p})$ then $J^{\ell_p}\subseteq J_T$.
\end{theorem}

\begin{proof}
Let $T\not\in\mathcal{SS}(d_{w,p})$. Then there exists an infinite-dimensional subspace $Y$ of $d_{w,p}$ such that $T|_Y$ is an isomorphism. By Remark~\ref{ell-p-wealth}, passing to a subspace, we may assume that $Y$ is complemented in $d_{w,p}$ and isomorphic to~$\ell_p$. Let $(x_n)$ be a basis of $Y$ equivalent to the unit vector basis of~$\ell_p$. Define $z_n=Tx_n$, then $(z_n)$ is also equivalent to the unit vector basis of~$\ell_p$. By Remark~\ref{ell-p-wealth}, $(z_n)$ has a subsequence $(z_{n_k})$ such that $[z_{n_k}]$ is complemented in $d_{w,p}$ and isomorphic to~$\ell_p$.

Denote $W=[x_{n_k}]$. Then $W$ and $T(W)$ are both complemented subspaces of $d_{w,p}$ isomorphic to~$\ell_p$. Let $P$ and $Q$ be projections onto $W$ and $T(W)$, respectively. Put $S=\big(T|_W\big)^{-1}$, $S\in L\big(T(W),d_{w,p}\big)$. Then it is easy to see that $P=(SQ)TP$. Since $SQ$ and $P$ are in $L(d_{w,p})$, we have $J_P\subseteq J_T$. By Proposition~\ref{ell-p-fact-characterization}, $J^{\ell_p}\subseteq J_T$.
\end{proof}

\begin{corollary} $\overline{J^{\ell_p}}\bigvee\mathcal{SS}(d_{w,p})$ is the only immediate successor of $\mathcal{SS}(d_{w,p})$ and
$\overline{J^{\ell_p}}$ is an immediate successor of $\overline{J^{\ell_p}}\wedge\mathcal{SS}(d_{w,p})$.
\end{corollary}

Now we will investigate the ideal of finitely strictly singular operators on~$d_{w,p}$. To prove the main statement (Theorem~\ref{FSS}), we will need the following lemma due to Milman~\cite{Mil70} (see also a thorough discussion in~\cite{SSTT}). This lemma will be used more than once in the paper.

\begin{lemma}[{\cite{Mil70}}]\label{milman}
If $F$ is a $k$-dimensional subspace of $c_0$ then there exists a vector $x\in F$ such that $x$ attains its sup-norm at at least $k$ coordinates (that is, $x^*$ starts with a constant block of length $k$).
\end{lemma}

We will also use the following simple lemma.

\begin{lemma}\label{norm-entry-correspondence}
Let $s_n=\sum_{i=1}^nw_i$ ($n\in\mathbb N$) where $w=(w_i)$ is the sequence of weights for~$d_{w,p}$. If $x\in d_{w,p}$, $y=x^*$, and $N\in\mathbb N$ then $0\le y_N\le\frac{\norm{x}}{s^{1/p}_N}$.
\end{lemma}
\begin{proof}
$\norm{x}^p=\norm{y}^p=\sum_{i=1}^\infty y^p_iw_i\ge y_N^p\sum_{i=1}^Nw_i=y_N^ps_N$.
\end{proof}

\begin{theorem}\label{FSS}
Let $X$ and $Y$ be subspaces of~$d_{w,p}$. Then $\mathcal{FSS}(X,Y)=\mathcal{SS}(X,Y)$. In particular, $\mathcal{FSS}(\ell_p,d_{w,p})=\mathcal{SS}(\ell_p,d_{w,p})$ and $\mathcal{FSS}(d_{w,p})=\mathcal{SS}(d_{w,p})$.
\end{theorem}

\begin{proof}
Let $T\in L(X,Y)$. Suppose that $T$ is not finitely strictly singular. We will show that it is not strictly singular. Since $T$ is not finitely strictly singular, there exists a constant $c>0$ and a sequence $F_n$ of subspaces of $X$ with $\dim F_n\ge n$ such that for each $n$ and for all $x\in F_n$ we have $\norm{Tx}\ge c\norm{x}$.

Fix a sequence $(\varepsilon_k)$ in $\mathbb R$ such that $1>\varepsilon_k\downarrow 0$. We will inductively construct a sequence $(x_k)$ in $X$ and two strictly increasing sequences $(n_k),(m_k)$ in $\mathbb N$ such that:
\begin{enumerate}
\item\label{i-seminorm} $(x_k)$ and $(Tx_k)$ are seminormalized; we will denote $Tx_k$ by $u_k$;
\item\label{ii-almost-disj} for all $k\in\mathbb N$, $\supp x_k\subseteq[n_k,\infty)$ and $\supp u_k\subseteq[m_k,\infty)$;
\item\label{u'-long} if $k\ge 2$ then $\norm{x_{k-1}|_{[n_k,\infty)}}<\varepsilon_k$, $\norm{u_{k-1}|_{[m_k,\infty)}}<\varepsilon_k$, and all the coordinates of $u_{k-1}$ where the sup-norm is attained are less than $m_k$;
\item\label{iii-long} for each $k\in\mathbb N$, the vector $u_k^*$ begins with a constant block of length at least~$k$.
\end{enumerate}
That is, $(x_n)$ and $(u_n)$ are two almost disjoint sequences and $u_n$'s have long ``flat'' sections.

Take $x_1$ to be any nonzero vector in $F_1$ and put $n_1=m_1=1$. Suppose we have already constructed $x_1,\dots,x_{k-1}$, $n_1,\dots,n_{k-1}$, and $m_1,\dots,m_{k-1}$ such that the conditions \eqref{i-seminorm}--\eqref{iii-long} are satisfied. Choose $n_k\in\mathbb N$ and $m_k\in\mathbb N$ such that $n_k>n_{k-1}$, $m_k>m_{k-1}$ and the condition~\eqref{u'-long} is satisfied.

Consider the space 
$$
V=\bigl\{y=(y_i)\in F_{n_k+m_k+k}\mid y_i=0\mbox{ for }i<n_k\bigr\}\subseteq F_{n_k+m_k+k}.
$$ 
It follows from $\dim F_{n_k+m_k+k}\ge n_k+m_k+k$ that $\dim V\ge m_k+k$. Since $V\subseteq F_{n_k+m_k+k}$, $\norm{Ty}\ge c\norm{y}$ for all $y\in V$. In particular, $\dim(TV)\ge m_k+k$. Define 
$$
Z=\bigl\{z=(z_i)\in TV\mid z_i=0\mbox{ for }i<m_k\bigr\}.
$$ 
It follows that $\dim Z\ge k$.

Clearly, $\supp y\subseteq [n_k,\infty)$ for all $y\in V$ and $\supp z\subseteq [m_k,\infty)$ for all $z\in Z$. By Lemma~\ref{milman}, we can choose $u_k\in Z$ such that $u_k$ is normalized and $u_k^*$ starts with a constant block of length~$k$. Put $x_k=(T|_V)^{-1}(u_k)\in Y$. Since $x_k\in V$ and $\norm{u_k}=1$, it follows that $\frac{1}{\norm{T}}\le\norm{x_k}\le\frac{1}{c}$, so the conditions \eqref{i-seminorm}--\eqref{iii-long} are satisfied for~$(x_k)$.

For each $k\in\mathbb N$, let $x_k'=x_k|_{[n_{k},n_{k+1})}$ and $u_k'=u_k|_{[m_{k},m_{k+1})}$. Passing to tails of sequences, if necessary, we may assume that both $(x_k')$ and $(u_k')$ are seminormalized block sequences of~$(e_n)$.

Since the non-increasing rearrangement of each $u'_{k}$ starts with a constant block of length $k$ by~\eqref{u'-long}, the coefficients in $u'_{k}$ converge to zero by Lemma~\ref{norm-entry-correspondence}. Therefore, passing to a subsequence, we may assume by Remark~\ref{ell-p-wealth} that $(u'_{k})$ is equivalent to the unit vector basis $(f_n)$ of~$\ell_p$. Using Theorem~\ref{small-perturbations} and passing to a further subsequence, we may also assume that $(x_{k})\sim(x'_{k})$ and $(u_{k})\sim(u'_{k})$.

By Proposition~\ref{domination}, the sequence $(x'_{k})$ is dominated by $(f_n)$. Notice that the condition $u_k=Tx_k$ implies $(x_k)\succeq(u_k)$. Therefore, we get the following chain of dominations and equivalences of basic sequences:
$$
(f_n)\succeq(x'_{k})\sim(x_{k})\succeq(u_{k})\sim(u'_{k})\sim(f_n).
$$
It follows that all the dominations in this chain are, actually, equivalences. In particular, $(x_{k})\sim(u_{k})$. Thus, $T$ is an isomorphism on the space $[x_{k}]$, hence $T$ is not strictly singular.
\end{proof}

Recall that an operator $T$ on a Banach space $X$ is weakly compact if the image of the unit ball of $X$ under $T$ is relatively weakly compact. Alternatively, $T$ is weakly compact if and only if for every bounded sequence $(x_n)$ in $X$ there exists a subsequence $(x_{n_k})$ of $(x_n)$ such that $(Tx_{n_k})$ is weakly convergent.

If $1<p<\infty$ then $d_{w,p}$ is reflexive, and, hence, every operator in $L(d_{w,p})$ is weakly compact. In case $p=1$ we have the following.

\begin{theorem}\label{weakly-compact}
Let $T\in L(d_{w,1})$. Then $T$ is weakly compact if and only if $T$ is strictly singular.
\end{theorem}
\begin{proof}
Suppose that $T$ is strictly singular. We will show that $T$ is weakly compact.

Let $(x_n)$ be a bounded sequence in~$X$. By Rosenthal's $\ell_1$-theorem, there is a subsequence $(x_{n_k})$ of $(x_n)$ such that $(x_{n_k})$ is either equivalent to the unit vector basis $(f_n)$ of $\ell_1$ or is weakly Cauchy. In the latter case, $(Tx_{n_k})$ is also weakly Cauchy. If $(x_{n_k})\sim(f_n)$ then, since $T$ is strictly singular, $(Tx_{n_k})$ cannot have subsequences equivalent to $(f_n)$. Hence, using Rosenthal's theorem one more time and passing to a further subsequence, we may assume that, again, $(Tx_{n_k})$ is weakly Cauchy. Since $d_{w,1}$ is weakly sequentially complete, the sequence $(Tx_{n_k})$ is weakly convergent. It follows that $T$ is weakly compact.


Conversely, let $J$ be the closed ideal of weakly compact operators in $L(d_{w,1})$. By the first part of the proof, $J$ is a successor of $\mathcal{SS}(d_{w,1})$. Suppose that $J\ne\mathcal{SS}(d_{w,1})$. By Theorem~\ref{J-P-V-SS}, $J^{\ell_1}\subseteq J$. This, however, is a contradiction since a projection onto a copy of $\ell_1$ (which belongs to $J^{\ell_1}$ by Proposition~\ref{ell-p-fact-characterization}) is not weakly compact.
\end{proof}

\section{Operators factorable through the formal identity}\label{form-dentity}

The operator $j\colon\ell_p\to d_{w,p}$ defined by $j(e_n)=f_n$ is called \term{the formal identity operator from $\ell_p$ to $d_{w,p}$}.
It follows immediately from the definition of the norm in $d_{w,p}$ that  $\norm{j}=1$.

We will denote by the symbol $J^j$ the set of all operators $T\in L(d_{w,p})$ which can be factored as $T=AjB$ where $A\in L(d_{w,p})$ and $B\in L(d_{w,p},\ell_p)$.

\begin{proposition}\label{J-j-ideal}
$J^j$ is an ideal in $L(d_{w,p})$.
\end{proposition}

\begin{proof}
It is clear from the definition that the set $J^j$ is closed under both right and left multiplication by operators from $L(d_{w,p})$. We have to show that if $T_1$ and $T_2$ are in $J^j$ then $T_1+T_2$ is in~$J^j$, as well. 

Write $T_1=A_1jB_1$, $T_2=A_2jB_2$ with $A_1,A_2\in L(d_{w,p})$ and $B_1,B_2\in L(d_{w,p},\ell_p)$. Let $A\in L(d_{w,p}\oplus d_{w,p},d_{w,p})$ and $B\in L(d_{w,p},\ell_p\oplus\ell_p)$ be defined by 
$$
A(x_1,x_2)=A_1x_1+A_2x_2\quad\mbox{and}\quad Bx=(B_1x,B_2x).
$$
Define also $U\colon\ell_p\to\ell_p\oplus\ell_p$ and $V\colon d_{w,p}\to d_{w,p}\oplus d_{w,p}$ by 
$$
U\big((x_n)\big)=\big((x_{2n-1}),(x_{2n})\big),\quad\mbox{and}\quad
V\big((x_n)\big)=\big((x_{2n-1}),(x_{2n})\big).
$$
Since the bases of $\ell_p$ and $d_{w,p}$ are both unconditional, $U$ and $V$ are bounded.

Now observe that for each $x=(x_n)\in d_{w,p}$ we can write 
\begin{eqnarray*}
AVjU^{-1}Bx=AVjU^{-1}(B_1x,B_2x)=&\\
A(jB_1x,jB_2x)=A_1jB_1x&\!\!\!\!+A_2jB_2x=T_1x+T_2x.
\end{eqnarray*}
This shows that $T_1+T_2=AVjU^{-1}B$ with $AV\in L(d_{w,p})$ and $U^{-1}B\in L(d_{w,p},\ell_p)$, hence $T_1+T_2\in J^j$.
\end{proof}

As we already mentioned before, the space $d_{w,p}$ contains many complemented copies of~$\ell_p$. Consider the operator $jUP\in L(d_{w,p})$ where $P$ is a projection onto any subspace $Y$ isomorphic to $\ell_p$ and $U\colon Y\to\ell_p$ is an isomorphism onto. It turns out that the ideal generated by any such operator does not depend on the choice of $Y$ and, in fact, coincides with~$J^j$.

\begin{proposition}\label{jUP}
Let $Y$ be a complemented subspace of $d_{w,p}$ isomorphic to~$\ell_p$, $P\in L(d_{w,p})$ be a projection with range~$Y$, and $U\colon Y\to\ell_p$ be an isomorphism onto. If $T=jUP$ then $J_T=J^j$.
\end{proposition}

\begin{proof}
Clearly, $J_T\subseteq J^j$. Let $S\in J^j$. Then $S=AjB$ where $A\in L(d_{w,p})$ and $B\in L(d_{w,p},\ell_p)$. It follows that 
$$
S=AjB=Aj(UPU^{-1})B=AT(U^{-1}B)\in J_T.
$$
\end{proof}

The next goal is to show that the ideal $\overline{J^j}$ ``sits'' between $\mathcal K(X)$ and $\mathcal{SS}(X)\wedge\overline{J^{\ell_p}}$.

\begin{theorem}\label{j-FSS}
The formal identity operator $j\colon\ell_p\to d_{w,p}$ is finitely strictly singular.
\end{theorem}

\begin{proof}
Let $\varepsilon>0$ be arbitrary. Take $n\in\mathbb N$ such that $\frac{1}{n}\sum_{i=1}^nw_i<\varepsilon$; such $n$ exists by $w_n\to 0$. Since $(w_n)$ is also a decreasing sequence, it follows that $w_i<\varepsilon$ for all $i\ge n$.

Let $Y\subseteq\ell_p$ be a subspace with $\dim Y\ge n$. By Lemma~\ref{milman}, there exists a vector $x\in Y$ such that $\norm{x}_{\ell_p}=1$ and $x$ attains its sup-norm at at least $n$ coordinates. Denote $\delta=\norm{x}_{\sup}>0$. Then $\norm{x}_{\ell_p}\ge n^{1/p}\delta$, so $\delta\le n^{-1/p}$. 

Observe that the non-increasing rearrangement $x^*$ of $x$ satisfies the condition that $x^*_i=\delta$ for all $1\le 1\le n$. Therefore
$$
\norm{jx}^p_{d_{w,p}}=\sum_{i=1}^\infty x^{*p}_iw_i\le\delta^p\sum_{i=1}^n w_i+\varepsilon\sum_{i=n+1}^\infty x^{*p}_i\le\delta^pn\varepsilon+\varepsilon\norm{x}^p_{\ell_p}\le 2\varepsilon.
$$
Hence $\norm{jx}_{d_{w,p}}\le(2\varepsilon)^{1/p}$.
\end{proof}

\begin{corollary}\label{J-j-compact}
The following inclusions hold: $\mathcal K(d_{w,p})\subsetneq\overline{J^j}$ and $J^j\subseteq \mathcal{SS}(d_{w,p})\wedge J^{\ell_p}$.
\end{corollary}

\begin{proof}
Let $Y$, $P$, and $U$ be as in Proposition~\ref{jUP}. Then $jUP\in J^j$. If $x_n=U^{-1}f_n\in d_{w,p}$ then $(x_n)$ is seminormalized and $jUPx_n=e_n$. Hence the sequence $(jUPx_n)$ has no convergent subsequences, so that $jUP$ is not compact.

The inclusion $J^j\subseteq \mathcal{SS}(d_{w,p})\wedge J^{\ell_p}$ is obvious since $j$ is strictly singular.
\end{proof}

\begin{conjecture}\label{K-J-imm}
  The ideal $\overline{J^j}$ is the only immediate successor of  $\mathcal K(d_{w,p})$.
\end{conjecture}

In \cite{ACL} and \cite{CL} (see also \cite{LT77}), conditions on the weights $w=(w_n)$ are given under which $d_{w,p}$ has exactly two non-equivalent symmetric basic sequences. We will show that the conjecture holds true in this case.

\begin{lemma}\label{T-x-n-w-null}
If $T\in \mathcal{SS}(d_{w,p})\setminus\mathcal K(d_{w,p})$ then there exists a seminormalized basic sequence $(x_n)$ in $d_{w,p}$ such that $(f_n)\succeq(x_n)$ and $(Tx_n)$ is weakly null and seminormalized.
\end{lemma}
\begin{proof}
Let $(z_n)$ be a bounded sequence in $d_{w,p}$ such that $(Tz_n)$ has no convergent subsequences. Then $(z_n)$ has no convergent subsequences either. Applying Rosenthal's $\ell_1$-theorem and passing to a subsequence, we may assume that $(z_n)$ is either equivalent to the unit vector basis of $\ell_1$ or is weakly Cauchy.

{\it Case: $(z_n)$ is equivalent to the unit vector basis of $\ell_1$}. Since a reflexive space cannot contain a copy of~$\ell_1$, we conclude that $p=1$, so $(z_n)\sim(f_n)$. Again, by Rosenthal's theorem, $(Tz_n)$ has a subsequence which is either equivalent to $(f_n)$ or is weakly Cauchy. If $(Tz_{n_k})\sim(f_n)$ then $T$ is an isomorphism on the space $[z_{n_k}]$, contrary to the assumption that $T\in\mathcal{SS}(d_{w,p})$. Therefore, $(Tz_{n_k})$ is weakly Cauchy. Put $x_{k}=z_{n_{2k}}-z_{n_{2k-1}}$. Then $(x_k)$ is basic and $(Tx_k)$ is weakly null. Passing to a further subsequence of $(z_{n_k})$ we may assume that $(Tx_k)$ is seminormalized. Also, $(x_k)$ is still equivalent to $(f_n)$, hence is dominated by~$(f_n)$.

{\it Case: $(z_n)$ is weakly Cauchy}. Clearly, $(Tz_n)$ is also weakly Cauchy. Consider the sequence $(u_n)$ in $d_{w,p}$ defined by $u_n=z_{2n}-z_{2n-1}$. Then both $(u_n)$ and $(Tu_n)$ are weakly null. Passing to a subsequence of $(z_n)$, we may assume that $(Tu_n)$ and, hence, $(u_n)$ are seminormalized. Applying Theorem~\ref{bessaga-pelczynski}, we get a subsequence $(u_{n_k})$ of $(u_n)$ which is basic and equivalent to a block sequence $(v_n)$ of $(e_n)$. Denote $x_k=u_{n_k}$. By Proposition~\ref{domination}, $(f_n)$ dominates $(v_n)$ and, hence, $(x_k)$.
\end{proof}

\begin{theorem}\label{J-j-direct}
If $d_{w,p}$ has exactly two non-equivalent symmetric basic sequences, then $\overline{J^j}$ is the only immediate successor of $\mathcal K(d_{w,p})$.
\end{theorem}

\begin{proof}
Let $T$ be a non-compact operator on~$d_{w,p}$. It suffices to show that $J^j\subseteq J_T$. We may assume that $T$ is strictly singular because, otherwise, we have $J^j\subseteq J^{\ell_p}\subseteq J_T$ by Theorem~\ref{J-P-V-SS}.

Let $(x_n)$ be a sequence as in Lemma~\ref{T-x-n-w-null}. Passing to a subsequence and using Theorem~\ref{bessaga-pelczynski}, we may assume that $(Tx_n)$ is basic and equivalent to a block sequence $(h_n)$ of $(e_n)$ such that $Tx_n-h_n\to 0$. We claim that $(h_n)$ has no subsequences equivalent to~$(f_n)$. Indeed, otherwise, for such a subsequence $(h_{n_k})$ of $(h_n)$, we would have $(f_n)\sim(f_{n_k})\succeq(x_{n_k})\succeq(Tx_{n_k})\sim(h_{n_k})\sim(f_n)$, so $(x_{n_k})\sim(Tx_{n_k})$, contrary to $T\in\mathcal{SS}(d_{w,p})$. By \cite[Theorem~19]{CL}, $(h_n)$ has a subsequence which spans a complemented subspace in $d_{w,p}$ and is equivalent to $(e_n)$. Therefore, by Theorem~\ref{small-perturbations}, we may assume (by passing to a further subsequence) that $(Tx_n)\sim(e_n)$ and $[Tx_{n}]$ is complemented in~$d_{w,p}$.

We have proved that there exists a sequence $(x_n)$ in $d_{w,p}$ such that $[Tx_n]$ is complemented in $d_{w,p}$ and
$$
(f_n)\succeq (x_n)\succeq (Tx_n)\sim(e_n).
$$
Let $A\in L(\ell_p,d_{w,p})$ and $B\in L([Tx_n],d_{w,p})$ be defined
by $Af_n=x_n$ and $B(Tx_n)=e_n$. Let $Q\in L(d_{w,p})$ be a projection
onto $[Tx_n]$. Then for all $n\in\mathbb N$, we obtain:
$BQTAf_n=BQTx_n=BTx_n=e_n$. It follows that $BQTA=j$, so that
$J^j\subseteq J_T$.
\end{proof}

In order to prove Conjecture~\ref{K-J-imm} without additional conditions on~$w$, it suffices to show that if $T\in\overline{J^j}\setminus\mathcal K(d_{w,p})$ then $J^j\subseteq\overline{J_T}$. We will prove a weaker statement: if $T\in J^j\setminus\mathcal K(d_{w,p})$ then $J^j\subseteq J_T$.

Recall (see \cite[p.148]{ACL}) that if $x=(a_n)\in d_{w,p}$ then a block sequence $(y_n)$ of $(e_n)$ is called a \term{block of type I generated by} $x$ if it is of the form $y_n=\sum_{i=p_n+1}^{p_{n+1}}a_{i-p_n}e_i$ for all~$n$. A set $A\subseteq d_{w,p}$ will be said to be \term{almost lengthwise bounded} if for each $\varepsilon>0$ there exists $N\in\mathbb N$ such that $\norm{x^*|_{[N,\infty)}}<\varepsilon$ for all $x\in A$. We will usually use it in the case when $A=\{x_n\}$ for some sequence $(x_n)$ in~$d_{w,p}$. We need the following result, which is a slight extension of \cite[Theorem~3]{ACL}. We include the proof for completeness.

\begin{theorem}\label{almost-bdd-thm}
Let $(x_n)$ be a seminormalized block sequence of $(e_n)$ in~$d_{w,p}$. 
\begin{enumerate}
\item\label{not-alb} If $(x_n)$ is not almost lengthwise bounded then there exists a subsequence $(x_{n_k})$ such that $(x_{n_k})\sim(f_n)$.
\item\label{alb} If $(x_n)$ is almost lengthwise bounded, then there exists a subsequence $(x_{n_k})$ equivalent to a block of type I generated by a vector $u=\sum_{i=1}^{\infty}b_ie_i\in d_{w,p}$ with $b_i\downarrow 0$. Moreover, if the sequence $(x_n)$ is bounded in $\ell_p$ then $u$ is in~$\ell_p$.\footnote{As a sequence space, $\ell_p$ is a subset of~$d_{w,p}$. That is, we can identify $\ell_p$ with $\Range j$. More precisely, we claim here that if $(j^{-1}x_n)$ is bounded in $\ell_p$ then $u$ is in $\Range j$. Being a block sequence of $(e_n)$, $(x_n)$ is contained in $\Range j$.}
\end{enumerate}
\end{theorem}

\begin{proof}
\eqref{not-alb}
Without loss of generality, $\norm{x_n}\le 1$ for all $n\in\mathbb N$. By the assumption, there exists $\varepsilon>0$ with the property that for each $k\in\mathbb N$, there is $n_k\in\mathbb N$ such that $\norm{x^*_{n_k}|_{(k,\infty)}}\ge\varepsilon$. Let $u_k$ be a restriction of  $x_{n_k}$ such that $u^*_k=x^*_{n_k}|_{[1,k]}$ and $v_k=x_{n_k}-u_k$.

Clearly, each nonzero entry of $u_k$ is greater than or equal to the greatest entry of~$v_k$. By Lemma~\ref{norm-entry-correspondence}, the $k$-th coordinate of $u^*_k$ is less than or equal to $\frac{1}{s_k^{1/p}}$ where $s_k=\sum_{i=1}^kw_i$. It follows that $(v_k)$ is a block sequence of $(e_n)$ such that $\varepsilon\le\norm{v_k}\le 1$ and absolute values of the entries of $v_k$ are all at most $\frac{1}{s_k^{1/p}}$. Since $\lim_ks_k=+\infty$ by the definition of~$d_{w,p}$, passing to a subsequence and using Remark~\ref{ell-p-wealth} we may assume that $(v_k)$ is equivalent to $(f_n)$. By Proposition~\ref{domination}, $(f_n)$ dominates $(x_{n_k})$. Using also Lemma~\ref{lattice-domination}, we obtain the following diagram:
$$
(f_n)\succeq(x_{n_k})\succeq(v_k)\sim(f_n).
$$
Hence $(x_{n_k})$ is equivalent to $(f_n)$.

\eqref{alb}
Suppose that $x_n=\sum_{i=p_n+1}^{p_{n+1}}a_ie_i$. Clearly, the sequence $(a_i)$ is  bounded. Without loss of generality, $a_{p_n+1}\ge\dots\ge a_{p_{n+1}}\ge 0$ for each~$n$. Put $y_n=x_n^*$. Using a standard diagonalization argument and passing to a subsequence, we may assume that $(y_n)$ converges coordinate-wise; put $b_i=\lim\limits_{n\to\infty}y_{n,i}$. It is easy to see that $b_i\ge b_{i+1}$ for all~$i$. Put $u=(b_i)$.

{\it Case: the sequence $(p_{n+1}-p_n)$ is bounded.} Passing to a subsequence, we may assume that $N:=p_{n_k+1}-p_{n_k}$ is a constant. Note that $\supp u\subseteq[1,N]$ and  $\supp y_{n_k}\subseteq[1,N]$ for all~$k$. Put $u_k=\sum_{i=p_{n_k}+1}^{p_{n_k+1}}b_{i-p_{n_k}}e_i$, then $u=u_k^*$ and $(u_k)$ as a block of type I generated by~$u$. By compactness, 
\begin{math}
\norm{x_{n_k}-u_k}
=\norm{y_{n_k}-u}\to 0.
\end{math}
Therefore, passing to a further subsequence, we have $(x_{n_k})\sim (u_k)$. Being a vector with finite support, $u$ belongs to~$\ell_p$.

{\it Case: the sequence $(p_{n+1}-p_n)$ is unbounded.} We will construct the required subsequence $(x_{n_k})$ and a sequence $(N_k)$ inductively. Put $n_1=N_1=1$ and let $k>1$. Suppose that $n_1,\dots,n_{k-1}$ and $N_1,\dots,N_{k-1}$ have already been selected. Since $(x_n)$ is almost lengthwise bounded, we can find $N_k>N_{k-1}$ such that $\norm{y_n|_{(N_k,\infty)}}<\frac{1}{k}$ for all~$n$. Put $v_k=u|_{[1,N_k]}$. Using coordinate-wise convergence, we can find $n_k>n_{k-1}$ such that $\norm{y_{n_k}|_{[1,N_k]}-v_k}_{\ell_p}<\frac{1}{k}$ and $p_{n_k}+N_k\le p_{n_k+1}$. Put $u_k=\sum_{i=p_{n_k}+1}^{p_{n_k}+N_k}b_{i-p_{n_k}}e_i$. Then $u_k^*=v_k$, so that
\begin{equation}
  \label{eq:lp-est}
     \norm{x_{n_k}|_{(p_{n_k},p_{n_k}+N_k]}-u_k}_{\ell_p}
    =\norm{y_{n_k}|_{[1,N_k]}-v_k}_{\ell_p}<\tfrac{1}{k}
\end{equation}
and
  \begin{displaymath}
    \norm{x_{n_k}|_{(p_{n_k}+N_k,p_{n_k+1}]}}
    =\norm{y_{n_k}|_{(N_k,\infty)}}<\tfrac{1}{k}.
  \end{displaymath}
It follows that $\norm{x_{n_k}-u_k}\to 0$. Passing to a subsequence, we get $(x_{n_k})\sim(u_k)$.

Next, we show that $u\in d_{w,p}$. Since $\norm{\cdot}\le\norm{\cdot}_{\ell_p}$, it follows from~\eqref{eq:lp-est} that
\begin{displaymath}
  \norm{v_k}=\norm{u_k}
  \le\norm{x_{n_k}|_{(p_{n_k},p_{n_k}+N_k]}}+\tfrac{1}{k}
  \le\norm{x_{n_k}}+\tfrac{1}{k}.
\end{displaymath}
Since $(x_n)$ is bounded, so is $(v_k)$. Since $\supp v_k=N_k\to\infty$, we have $u\in d_{w,p}$. For the ``moreover'' part, we argue in a similar way. By~\eqref{eq:lp-est}, we have 
\begin{displaymath}
  \norm{v_k}_{\ell_p}
  \le\norm{u_k}_{\ell_p}
  \le\norm{x_{n_k}|_{(p_{n_k},p_{n_k}+N_k]}}_{\ell_p}+\tfrac{1}{k}
  \le\norm{x_{n_k}}_{\ell_p}+\tfrac{1}{k}.
\end{displaymath}
Therefore, if $(x_n)$ is bounded in $\ell_p$ then so is $(v_k)$, hence $u\in\ell_p$. 
\end{proof}

\begin{lemma}\label{type-I-e-n}
Suppose that $(u_n)$ is a block of type I in $d_{w,p}$ generated by some $u=\sum_{i=1}^\infty b_ie_i$. If $b_i\downarrow 0$ and $u\in\ell_p$ then $(u_n)$ has a subsequence equivalent to $(e_n)$
\end{lemma}

\begin{proof}
  By Corollary~4 of~\cite{ACL}, we may assume that the basic sequence $(u_n)$ is symmetric. It suffices to show that $[u_n]$ is isomorphic to $d_{w,p}$ because all symmetric bases in $d_{w,p}$ are equivalent; see e.g., Theorem~4 of~\cite{ACL}. Without loss of generality, $\norm{u}=1$. Lemma~4 of~\cite{ACL} asserts that $[u_n]$ is isomorphic to $d_{w,p}$ iff $(s_n^{(u)})\sim (s_n)$, where $s_n=\sum_{i=1}^nw_i$, $s_n^{(u)}=\sum_{i=1}^\infty b_i^p(s_{ni}-s_{n(i-1)})$, and $(\alpha_n)\sim(\beta_n)$ means that there exist positive constants $A$ and $B$ such that $A\alpha_n\le \beta_n\le B\alpha_n$ for all~$n$. Let's verify that this condition is, indeed, satisfied. On one hand, taking only the first term in the definition of $s_n^{(u)}$, we get $s_n^{(u)}\ge b_1^ps_n$. On the other hand, it follows from $w_i\downarrow$ that $s_{ni}-s_{n(i-1)}\le s_n$ for every~$i$, hence $s_n^{(u)}\le\sum_{i=1}^\infty b_i^ps_n=\norm{u}_{\ell_p}^ps_n$.
\end{proof}

\begin{lemma}\label{jxn-equiv-en}
Let $(x_n)$ be a block sequence of $(f_n)$ in $\ell_p$ such that the sequences $(x_n)$ and $(jx_n)$ are seminormalized in $\ell_p$ and~$d_{w,p}$, respectively. Then there exists a subsequence $(x_{n_k})$ such that $(jx_{n_k})\sim(e_n)$.
\end{lemma}

\begin{proof}
Clearly, $(x_n)\sim(f_n)$. It follows that $(jx_n)\not\sim(f_n)$ because, otherwise, $j$ would be an isomorphism on $[x_n]$, which is impossible because $j$ is strictly singular by Theorem~\ref{j-FSS}. Applying Theorem~\ref{almost-bdd-thm} to $(jx_n)$ and passing to a subsequence, we may assume that $(jx_n)\sim(u_n)$, where $(u_n)$ is a block of type I generated by some $u=\sum_{i=1}^\infty b_ie_i$ such that $b_i\downarrow 0$ and $u\in\ell_p$. Applying Lemma~\ref{type-I-e-n} and passing to a subsequence, we get $(u_n)\sim(e_n)$.
\end{proof}

\begin{theorem}
If $T\in J^j\setminus\mathcal K(d_{w,p})$ then $J^j\subseteq J_T$.
\end{theorem}
\begin{proof}
Write $T=AjB$ where $B\colon d_{w,p}\to\ell_p$ and $A\colon d_{w,p}\to d_{w,p}$. Let $(x_n)$ be as in Lemma~\ref{T-x-n-w-null}. The sequence $(Bx_n)$ is bounded, hence we may assume by passing to a subsequence that it converges coordinate-wise. Since $(Tx_n)$ is weakly null and seminormalized, it has no convergent subsequences. It follows that, after passing to a subsequence of $(x_n)$, we may assume that $(Tz_n)$ is seminormalized, where $z_n=x_{2n}-x_{2n-1}$. In particular, $(z_n)$, $(Bz_n)$, and $(jBz_n)$ are seminormalized. Also, $(Bz_n)$ converges to zero coordinate-wise. Using Theorem~\ref{bessaga-pelczynski} and passing to a further subsequence, we may assume that $(Bz_n)$  is equivalent to a block sequence $(u_n)$ of $(f_n)$ and $Bz_n-u_n\to 0$. It follows from $(f_n)\succeq(x_n)$ that $(f_n)\succeq(z_n)\succeq(Bz_n)\sim(u_n)\sim(f_n)$. In particular, $(z_n)\sim(f_n)$.

Since $Bz_n-u_n\to 0$ and $(jBz_n)$ is seminormalized, we may assume that the sequence $(ju_n)$ is seminormalized. By Lemma~\ref{jxn-equiv-en}, passing to a further subsequence, we may assume that $(ju_n)$ and, hence, $(jBz_n)$ are equivalent to $(e_n)$.

Passing to a subsequence and using Theorem~\ref{bessaga-pelczynski}, we may assume that $(Tz_n)$ is equivalent to a block sequence $(v_n)$ of $(e_n)$ such that $Tz_n-v_n\to 0$. Since $T\in\mathcal{SS}(d_{w,p})$, no subsequence of $(Tz_n)$ and, therefore, of $(v_n)$ is equivalent to $(f_n)$. By Proposition~\ref{domination}, $(v_n)\succeq(e_n)$. It follows from $(jBz_n)\sim(e_n)$ that $(e_n)\succeq(Tz_n)$, hence $(Tz_n)\sim(e_n)\sim(v_n)$. 

Write $v_n=\sum_{i=p_n+1}^{p_{n+1}}a_ne_n$. By Remark~\ref{ell-p-wealth}, $a_n\not\to 0$. Hence, passing to a subsequence and using \cite[Remark~9]{CL}, we may assume that $[v_n]$ is complemented. By Theorem~\ref{bessaga-pelczynski}, we may assume that $[Tz_n]$ is complemented. Let $P\in L(d_{w,p})$ be a projection onto $[Tz_n]$ and $U\in L(\ell_p,d_{w,p})$ and $V\in L([Tz_n],d_{w,p})$ be defined by $Uf_n=z_n$ and $VTz_n=e_n$. Then we can write $j=VPTU$. Therefore $J^j\subseteq J_T$.
\end{proof}

\section{$d_{w,p}$-strictly singular operators}\label{S_E-section}

The ideals in $L(d_{w,p})$ we have obtained so far can be arranged into the following diagram.
$$
\xymatrix @R=0pt@C=20pt{ & & &\mathcal{SS}\ar@{=>}[dr]& & \\\{0\}\ar@{=>}[r]& \mathcal K\subsetneq\overline{J^j}\ar[r]&\overline{J^{\ell_p}}\wedge\mathcal{SS}\ar@{:>}[dr]\ar[ur]& &\overline{J^{\ell_p}}\vee\mathcal{SS}\ar[r]&L(d_{w,p})\\ & & &\overline{J^{\ell_p}}\ar[ur]& & }
$$
(see the Introduction for the notations). In this section, we will characterize the greatest ideal in the algebra $L(d_{w,p})$, that is, a proper ideal in $L(d_{w,p})$ that contains all other proper ideals in $L(d_{w,p})$.

If $X$ and $Y$ are two Banach spaces, then an operator $T\in L(X)$ is called \term{$Y$-strictly singular} if for any subspace $Z$ of $X$ isomorphic to~$Y$, the restriction $T|_Z$ is not an isomorphism.
The set of all $Y$-strictly singular operators in $L(d_{w,p})$ will be denoted by $\mathcal{SS}_Y$.

According to this notation, the symbol $\mathcal{SS}_{d_{w,p}}$ stands for the set of all $d_{w,p}$-strictly singular operators in $L(d_{w,p})$ (not to be confused with $\mathcal{SS}(d_{w,p})$). 

\begin{lemma}\label{Tx_n-converge}
Suppose that $T\in\mathcal{SS}_{d_{w,p}}$ and $(x_n)$ is a basic sequence in $d_{w,p}$ equivalent to the unit vector basis $(e_n)$. Then $Tx_n\to 0$.
\end{lemma}
\begin{proof}
Suppose, by way of contradiction, that $Tx_n\not\to 0$. Then there is a subsequence $(x_{n_k})$ such that $(Tx_{n_k})$ is seminormalized. Since $(x_n)$ is weakly null (Remark~\ref{e-n-w-null-remark}), we may assume by using Theorem~\ref{bessaga-pelczynski} and passing to a further subsequence that $(Tx_{n_k})$ is a basic sequence equivalent to a block sequence $(z_k)$ of~$(e_n)$.

By Proposition~\ref{domination}, either $(z_k)$ has a subsequence equivalent to $(f_n)$ or $(z_k)\succeq(e_n)$. Since $(Tx_{n_k})$ cannot have subsequences equivalent to $(f_n)$ (this would contradict boundedness of $T$), the former is impossible. Therefore $(z_k)\succeq(e_n)$. We obtain the following diagram:
$$
(e_n)\sim(x_{n_k})\succeq(Tx_{n_k})\sim(z_k)\succeq(e_n).
$$
Therefore $T|_{[x_{n_k}]}$ is an isomorphism. This contradicts $T$ being in $\mathcal{SS}_{d_{w,p}}$.
\end{proof}

\begin{corollary}\label{S_E-compact-characterize}
Let $T\in\mathcal{SS}_{d_{w,p}}$. If $Y\subseteq d_{w,p}$ is a subspace isomorphic to $d_{w,p}$ then there is a subspace $Z\subseteq Y$ such that $Z$ is isomorphic to $d_{w,p}$ and $T|_Z$ is compact.
\end{corollary}
\begin{proof}
Let $(x_n)$ be a basis of $Y$ equivalent to $(e_n)$. By Lemma~\ref{Tx_n-converge}, $Tx_n\to 0$. There is a subsequence $(x_{n_k})$ of $(x_n)$ such that $\sum_{k=1}^\infty\frac{\norm{Tx_{n_k}}}{\norm{x_{n_k}}}$ is convergent. Let $Z=[x_{n_k}]$. It follows that $Z$ is isomorphic to $d_{w,p}$ and $T|_Z$ is compact (see, e.g., \cite[Lemma 5.4.10]{CPY}).
\end{proof}

\begin{theorem}\label{SS-E-ideal}
The set $\mathcal{SS}_{d_{w,p}}$ of all $d_{w,p}$-strictly singular operators in $L(d_{w,p})$ is the greatest proper ideal in the algebra $L(d_{w,p})$. In particular, $\mathcal{SS}_{d_{w,p}}$ is closed.
\end{theorem}

\begin{proof}
First, let us show that $\mathcal{SS}_{d_{w,p}}$ is an ideal. Let $T\in\mathcal{SS}_{d_{w,p}}$. If $A\in L(d_{w,p})$ then, trivially, $AT\in\mathcal{SS}_{d_{w,p}}$. If $TA\not\in\mathcal{SS}_{d_{w,p}}$ then there exists a subspace $Y$ of $d_{w,p}$ such that $Y$ and $TA(Y)$ are both isomorphic to~$d_{w,p}$. Then $A|_Y$ is bounded below, hence $A(Y)$ is isomorphic to~$d_{w,p}$. It follows that $T$ is an isomorphism on a copy of~$d_{w,p}$, contrary to $T\in\mathcal{SS}_{d_{w,p}}$. So, $\mathcal{SS}_{d_{w,p}}$ is closed under two-sided multiplication by bounded operators.

Let $T,S\in\mathcal{SS}_{d_{w,p}}$. We will show that $T+S\in\mathcal{SS}_{d_{w,p}}$. Let $Y$ be a subspace of $d_{w,p}$ isomorphic to~$d_{w,p}$. By Corolary~\ref{S_E-compact-characterize}, there exists a subspace $Z$ of $Y$ such that $Z$ is isomorphic to $d_{w,p}$ and $T|_Z$ is compact. Applying Corolary~\ref{S_E-compact-characterize} again, we can find a subspace $V$ of $Z$ such that $V$ is isomorphic to $d_{w,p}$ and $S|_V$ is compact. Therefore $(T+S)|_V$ is compact, so that $(T+S)|_Y$ is not an isomorphism. So, $\mathcal{SS}_{d_{w,p}}$ is an ideal.

Clearly, the identity operator $I$ does not belong to $\mathcal{SS}_{d_{w,p}}$, so $\mathcal{SS}_{d_{w,p}}$ is proper. Let us show that $\mathcal{SS}_{d_{w,p}}$ is the greatest ideal in $L(d_{w,p})$.

Let $T\not\in\mathcal{SS}_{d_{w,p}}$. Then there exists a subspace $Y$ of $d_{w,p}$ such that $Y$ and $T(Y)$ are isomorphic to~$d_{w,p}$. By \cite[Corollary~12]{CL}, there exists a complemented (in $d_{w,p}$) subspace $Z$ of $T(Y)$ such that $Z$ is isomorphic to~$d_{w,p}$. Let $P\in L(d_{w,p})$ be a projection onto~$Z$. Put $H=T^{-1}(Z)$. It follows that $H$ is isomorphic to~$d_{w,p}$. Let $U\colon d_{w,p}\to H$ and $V\colon Z\to d_{w,p}$ be surjective isomorphisms. Then $S\in L(d_{w,p})$ defined by $S=(VP)TU$ is an invertible operator. Clearly $S\in J_T$, hence $J_T=L(X)$.

The fact that $\mathcal{SS}_{d_{w,p}}$ is closed follows from \cite[Corollary~VII.2.4]{Conway}.
\end{proof}

The next theorem provides a convenient characterization of $d_{w,p}$-strictly singular operators.

\begin{lemma}\label{Te-n-to-0}
Let $T\in L(d_{w,p})$ be such that $Te_n\to 0$. Suppose that $(x_n)$ is a bounded block  sequence of $(e_n)$ in $d_{w,p}$ such that $(x_n)$ is almost lengthwise bounded. Then $Tx_n\to 0$.
\end{lemma}

\begin{proof}
Write $x_n=\sum_{i=p_n+1}^{p_{n+1}}a_{i}e_i$. Since $(x_n)$ is bounded, there is $C>0$ such that $\abs{a_{i}}\le C$ for all $i$ and $n\in\mathbb N$. Let $\varepsilon>0$. Find $N\in\mathbb N$ such that $\norm{x^*_n|_{[N,\infty)}}<\varepsilon$ for all $n\in\mathbb N$. Let $u_n$ be a restriction of  $x_n$ such that $u^*_n=x^*_{n}|_{[1,N)}$ and $v_n=x_{n}-u_n$. It is clear that $\norm{v_n}=\norm{x^*_n|_{[N,\infty)}}<\varepsilon$. Also, $\norm{Tu_n}\le NC\cdot\max_{p_n+1\le i\le p_{n+1}}\norm{Te_i}$.

Pick $M\in\mathbb N$ such that $\norm{Te_k}<\frac{\varepsilon}{N}$ for all $k\ge M$. Then 
$$
\norm{Tx_n}\le\norm{Tu_n}+\norm{Tv_n}\le NC\frac{\varepsilon}{N}+\varepsilon\norm{T}=\varepsilon\big(C+\norm{T}\big)
$$
for all $n$ such that $p_n>M$. It follows that $Tx_n\to 0$. 
\end{proof}

\begin{theorem}\label{SS-E-characterization}
An operator $T\in L(d_{w,p})$ is $d_{w,p}$-strictly singular if and only if $Te_n\to 0$.
\end{theorem}

\begin{proof}
Suppose that $Te_n\to 0$ but $T\notin\mathcal{SS}_{d_{w,p}}$. Then there exists a subspace $Y$ of $d_{w,p}$ such that $Y$ is isomorphic to $d_{w,p}$ and $T|_Y$ is an isomorphism. Let $(x_n)$ be a basis of $Y$ equivalent to $(e_n)$. By Remark~\ref{e-n-w-null-remark}, $x_n\overset{w}{\longrightarrow}0$. Using Theorem~\ref{bessaga-pelczynski} and passing to a subsequence, we may assume that $(x_n)$ is equivalent to a block sequence $(z_n)$ of $(e_n)$ such that $x_n-z_n\to 0$. Since $(z_n)$ is equivalent to $(e_n)$, it is almost lengthwise bounded by Theorem~\ref{almost-bdd-thm}. By Lemma~\ref{Te-n-to-0}, $Tz_n\to 0$. Since $x_n-z_n\to 0$, we obtain $Tx_n\to 0$. This is a contradiction since $(x_n)$ is seminormalized and $T|_{[x_n]}$ is an isomorphism.

The converse implication follows from Lemma~\ref{Tx_n-converge}.
\end{proof}

\begin{remark}
In Theorem~\ref{SS-E-ideal} we showed, in particular, that $\mathcal{SS}_{d_{w,p}}$ is closed under addition. Alternatively, we could have deduced this from Theorem~\ref{SS-E-characterization}.
\end{remark}

Recall that an operator $T$ on a Banach space $X$ is called \term{Dunford-Pettis} if for any sequence $(x_n)$ in~$X$, $x_n\overset{w}{\longrightarrow} 0$ implies $Tx_n\to 0$. If $1<p<\infty$ then the class of Dunford-Pettis operators on $d_{w,p}$ coincides with $\mathcal K(d_{w,p})$ because $d_{w,p}$ is reflexive. For the case $p=1$ we have the following result.

\begin{theorem}\label{dunford-pettis}
Let $T\in L(d_{w,1})$. Then $T$ is $d_{w,1}$-strictly singular if and only if $T$ is Dunford-Pettis.
\end{theorem}
\begin{proof}
If $T$ is Dunford-Pettis then then $T$ is $d_{w,1}$-strictly singular by Theorem~\ref{SS-E-characterization} because $(e_n)$ is weakly null.

Conversely, suppose that $T$ is $d_{w,1}$-strictly singular. Let $(x_n)$ be a weakly null sequence. Suppose that $(Tx_n)$ does not converge to zero. Then, passing to a subsequence, we may assume that $(x_n)$ is a seminormalized weakly null basic sequence equivalent to a block sequence $(u_n)$ of $(e_n)$ such that $x_n-u_n\to 0$. Clearly, $(u_n)$ is weakly null. In particular, $(u_n)$ has no subsequences equivalent to $(f_n)$. By Theorem~\ref{almost-bdd-thm}, $(u_n)$ is almost lengthwise bounded. Hence, by Lemma~\ref{Te-n-to-0}, $Tu_n\to 0$. It follows that $Tx_n\to 0$, contrary to the choice of $(x_n)$.
\end{proof}

\section{Strictly singular operators between $\ell_p$ and $d_{w,p}$.}\label{additional-results}

We do not know whether the ideals $\overline{J^j}$, $\mathcal{SS}\wedge\overline{J^{\ell_p}}$, and $\mathcal{SS}$ are distinct. In this section, we discuss some connections between these ideals.

\begin{conjecture}\label{conjecture}$\overline{J^j}=\mathcal{SS}\wedge\overline{J_{\ell_p}}$. In particular, every strictly singular operator in $L(d_{w,p})$ which factors through $\ell_p$ can be approximated by operators that factor through~$j$.
\end{conjecture}

The following statement is a refinement of Lemma~\ref{lattice-domination}. Recall that $d_{w,p}$ is a Banach lattice with respect to the coordinate-wise order.

\begin{lemma}\label{modify-lattice-domination}
Suppose that $(x_n)$ and $(y_n)$ are seminormalized sequences in $d_{w,p}$ such that $\abs{x_n}\ge\abs{y_n}$ for all $n\in\mathbb N$ and $x_n\to 0$ coordinate-wise. Then there exists an increasing sequence $(n_k)$ in $\mathbb N$ such that $(x_{n_k})$ and $(y_{n_k})$ are basic and $(x_{n_k})\succeq(y_{n_k})$.
\end{lemma}
\begin{proof}
Clearly, $y_n\to 0$ coordinate-wise. By Theorem~\ref{bessaga-pelczynski}, we can find a sequence $(n_k)$ and two block sequences $(u_k)$ and $(v_k)$ of $(e_n)$ such that $(x_{n_k})$ and $(y_{n_k})$ are basic, $(x_{n_k})\sim(u_k)$, $(y_{n_k})\sim(v_k)$, $x_{n_k}-u_k\to 0$, $y_{n_k}-v_k\to 0$, and for each $k\in\mathbb N$, the vector $u_k$ ($v_k$, respectively) is a restriction of $(x_{n_k})$ (of $(y_{n_k})$, respectively).

For each $k\in\mathbb N$, define $h_k\in d_{w,p}$ by putting its $i$-th coordinate to be equal to $h_k(i)=\sign\big(v_k(i)\big)\cdot \big(\abs{u_k(i)}\wedge \abs{v_k(i)}\big)$. Then $(h_k)$ is a block sequence of $(e_n)$ such that $\abs{h_k}\le\abs{u_k}$. A straightforward verification shows that $\abs{h_k-v_k}\le\abs{u_k-x_{n_k}}$.
It follows that $h_k-v_k\to 0$. By Theorem~\ref{small-perturbations}, passing to a subsequence, we may assume that $(h_k)$ is basic and $(h_k)\sim(v_k)$. By Lemma~\ref{lattice-domination}, $(u_k)\succeq(h_k)$. Hence $(x_{n_k})\succeq(y_{n_k})$.
\end{proof}

The next lemma is a version of Theorem~\ref{almost-bdd-thm} for the case $(x_n)$ is an arbitrary bounded sequence.

\begin{lemma}\label{almost-bdd-ext}
If the bounded sequence $(x_n)$ in $d_{w,p}$ is not almost lengthwise bounded, then there is a subsequence $(x_{n_k})$ such that  $(x_{n_{2k}}-x_{n_{2k-1}})$ is equivalent to the unit vector basis $(f_n)$ of~$\ell_p$.
\end{lemma}
\begin{proof}
We can assume without loss of generality that no subsequence of $(x_n)$ is equivalent to the unit vector basis of~$\ell_1$. Indeed, if $(x_{n_k})$ is equivalent to the unit vector basis of $\ell_1$ then $p=1$. It follows that $(x_{n_k})$ is equivalent to $(f_n)$ and hence $(x_{n_{2k}}-x_{n_{2k-1}})$ is equivalent to $(f_n)$, as well.

Without loss of generality, $\sup_{n}\norm{x_n}=1$. Since $(x_n)$ is not almost lengthwise bounded, there exists $c>0$ such that
\begin{equation}
\label{eq:not-alm-bdd}
\forall N\in\mathbb N\quad \exists n\in\mathbb N\quad \norm{x^*_n|_{[N,\infty)}}>c.
\end{equation}
Let $\frac{c}{4}>\varepsilon_k\downarrow 0$. We will inductively construct increasing sequences $(n_k)$ and $(N_k)$ in $\mathbb N$ and a sequence $(y_k)$ in $d_{w,p}$ such that the following conditions are satisfied for each $k$:
\begin{enumerate}
\item\label{enum-i} $\norm{x_{n_k}|_{[N_{k+1},\infty)}}<\varepsilon_k$;
\item\label{enum-ii} $y_k$ is supported on $[N_k,N_{k+1})$;
\item\label{enum-iii} $y_k$ is a restriction of $x_{n_k}$;
\item\label{enum-iv} $\norm{y_k}>\frac{c}{2}$;
\item\label{enum-v} $\norm{y_k}_\infty\le s^{-1/p}_{N_{k}}$ where $s_N$ is as in Lemma~\ref{norm-entry-correspondence}.
\end{enumerate}
For $k=1$, we put $N_1=1$, and define $n_1$ to be the first number $n$ such that $\norm{x_n}>c$; such an $n$ exists by~\eqref{eq:not-alm-bdd}. Pick $N_2\in\mathbb N$ such that $\norm{x_{n_1}|_{[N_2,\infty)}}<\varepsilon_1$. Put $y_1=x_{n_1}|_{[N_1,N_2)}$. It follows that $1\ge\norm{y_1}>c-\varepsilon_1>\frac{c}{2}$, and the coordinates of $y_1$ are all at most $1$ ($=s^{-1/p}_{1}$), hence all the conditions \eqref{enum-i}--\eqref{enum-v} are satisfied for $k=1$.

Suppose that appropriate sequences $(n_i)_{i=1}^k$, $(N_i)_{i=1}^{k+1}$, and $(y_i)_{i=1}^k$ have been constructed. Use~\eqref{eq:not-alm-bdd} to find $n_{k+1}$ such that $\norm{x^*_{n_{k+1}}|_{[2N_{k+1},\infty)}}>c$. Let $z$ be the vector obtained from $x_{n_{k+1}}$ by replacing its $N_{k+1}$ largest (in absolute value) entries with zeros. Then $\norm{z|_{[N_{k+1},\infty)}}\ge\norm{z^*|_{[N_{k+1},\infty)}}=\norm{x^*_{n_{k+1}}|_{[2N_{k+1},\infty)}}>c$. By Lemma~\ref{norm-entry-correspondence}, $\norm{z}_\infty\le s^{-1/p}_{N_{k+1}}$. Choose $N_{k+2}$ such that $\norm{x_{n_{k+1}}|_{[N_{k+2},\infty)}}<\varepsilon_{k+1}$. It follows that $\norm{z|_{[N_{k+2},\infty)}}<\varepsilon_{k+1}$. Put $y_{k+1}=z|_{[N_{k+1},N_{k+2})}$. Then $\norm{y_{k+1}}\ge c-\varepsilon_{k+1}>\frac{c}{2}$, and the inductive construction is complete.

The sequence $(y_k)$ constructed above is a seminormalized block sequence of $(e_n)$ such that the coordinates of $(y_k)$ converge to zero by condition~\eqref{enum-v}. Using Remark~\ref{ell-p-wealth} and passing to a subsequence, we may assume that $(y_k)$ is equivalent to the unit vector basis $(f_n)$ of~$\ell_p$.

Since $(x_n)$ contains no subsequences equivalent to the unit vector basis of~$\ell_1$, using the Rosenthal's $\ell_1$-theorem and passing to a further subsequence, we may assume that $(x_{n_k})$ is weakly Cauchy. For all $m>k\in\mathbb N$, we have: $\norm{x_{n_k}|_{[N_m,\infty)}}\le\norm{x_{n_k}|_{[N_{k+1},\infty)}}\le\varepsilon_k$. Therefore $\norm{x_{n_m}-x_{n_k}}\ge\norm{(x_{n_m}-x_{n_k})|_{[N_m,\infty)}}\ge\norm{x_{n_m}|_{[N_m,\infty)}}-\varepsilon_k\ge\norm{y_m}-\varepsilon_k\ge\frac{c}{2}-\varepsilon_k>\frac{c}{4}$. It follows that the sequence $(u_k)$ defined by $u_k=x_{n_{2k}}-x_{n_{2k-1}}$ is seminormalized and weakly null. Passing to a subsequence of $(x_{n_k})$, we may assume that $(u_k)$ is equivalent to a block sequence of $(e_n)$. By Proposition~\ref{domination}, $(f_n)\succeq(u_k)$.

Let $v_k=x_{n_{2k}}-\big(x_{n_{2k-1}}|_{[1,N_{2k})}\,\big)$. Then $\norm{u_k-v_k}=\norm{x_{n_{2k-1}}|_{[N_{2k},\infty)}}<\varepsilon_{2k-1}\to 0$. By Theorem~\ref{small-perturbations}, passing to a subsequence of $(x_{n_k})$, we may assume that $(v_k)$ is basic and $(v_k)\sim(u_k)$. Also, $(v_k)$ is weakly null. Note that $\abs{y_{2k}}\le\abs{v_k}$ for all $k\in\mathbb N$, since $\supp y_{2k}\subseteq [N_{2k},N_{2k+1})$, so that $y_{2k}$ is a restriction of~$v_k$. By Lemma~\ref{modify-lattice-domination}, passing to a subsequence, we may assume that $(v_k)\succeq(y_{2k})$. Therefore we obtain the following diagram:
$$
(f_k)\succeq(u_k)\sim(v_k)\succeq(y_{2k})\sim(f_{2k})\sim(f_n).
$$
It follows that $(u_k)$ is equivalent to $(f_k)$.
\end{proof}

\begin{corollary}\label{ell-p-almost-bdd}
If $T\in \mathcal{SS}(\ell_p,d_{w,p})$ then the sequence $(Tf_n)$ is almost lengthwise bounded.
\end{corollary}
\begin{proof}
Suppose that $(Tf_n)$ is not almost lengthwise bounded. By Lemma~\ref{almost-bdd-ext}, there is a subsequence $(f_{n_k})$ such that $(Tf_{n_{2k}}-Tf_{n_{2k-1}})$ is equivalent to $(f_n)$. It follows that $T|_{[f_{n_{2k}}-f_{n_{2k-1}}]}$ is an isomorphism.
\end{proof}

\begin{remark}
If we view $T$ as an infinite matrix, the vectors $(Tf_n)$ represent its columns.
\end{remark}

\begin{theorem}\label{ell-1-factor}
If $T\in L(\ell_1,d_{w,1})$ is such that the sequence $(Tf_n)$ is almost lengthwise bounded, then for any $\varepsilon>0$ there exists $S\in L(\ell_1)$ such that $\norm{T-jS}<\varepsilon$, where $j\in L(\ell_1,d_{w,1})$ is the formal identity operator.
\end{theorem}
\begin{proof}
Let $\varepsilon>0$ be fixed. Find $N\in\mathbb N$ such that $\norm{(Tf_n)^*|_{[N,\infty)}}<\varepsilon$ for all~$n$. Let $z_n\in d_{w,1}$ be the vector obtained from $Tf_n$ by keeping its largest $N$ coordinates and replacing the rest of the coordinates with zeros. 

Define $S\colon\ell_1\to d_{w,1}$ by $Sf_n=z_n$. Note that $\norm{T-S}=\sup_{n}\norm{(T-S)f_n}=\sup_{n}\norm{Tf_n-z_n}\le\varepsilon$; in particular, $S$ is bounded. 
Let $F=\Span\{e_1,\dots,e_N\}$. Since $\dim F<\infty$, there exists $C>0$ such that 
$$
\frac{1}{C}\norm{x}_{\ell_1}\le\norm{x}_{d_{w,1}}\le C\norm{x}_{\ell_1}
$$ 
for all $x\in F$. Observe that for each $n\in\mathbb N$, the non-increasing rearrangement $(Sf_n)^*$ is in~$F$. Therefore, for all $n\in\mathbb N$, we have
$$
\norm{Sf_n}_{\ell_1}=\norm{(Sf_n)^*}_{\ell_1}\le C\norm{(Sf_n)^*}_{d_{w,1}}=C\norm{Sf_n}_{d_{w,1}}\le C\norm{S}.
$$
It follows that the operator $\widetilde S\colon\ell_1\to\ell_1$ defined by $\widetilde Sf_n=Sf_n$ belongs to $L(\ell_1)$. Obviously, $S=j\widetilde S$. So, $\norm{T-j\widetilde S}<\varepsilon$.
\end{proof}

The next corollary follows immediately from Theorem~\ref{ell-1-factor} and Corollary~\ref{ell-p-almost-bdd}. This corollary can be considered as a support for Conjecture~\ref{conjecture}.

\begin{corollary}\label{ell-1-ss-factor}
$\mathcal{SS}(\ell_1,d_{w,1})$ is contained in the closure of $\{jS\mid S\in L(\ell_1,d_{w,1})\}$.
\end{corollary}

\begin{question} Does Corollary~\ref{ell-1-ss-factor} remain valid for $p>1$?
\end{question}

The following fact is standard, we include its proof for convenience of the reader.

\begin{proposition}\label{SS-X-ell-1}
If $X$ is a Banach space then $\mathcal{SS}(X,\ell_1)=\mathcal K(X,\ell_1)$.
\end{proposition}
\begin{proof}
Let $T\not\in\mathcal{K}(X,\ell_1)$. Pick a bounded sequence $(x_n)$ in $X$ such that $(Tx_n)$ has no convergent subsequences. By Schur's theorem, $(Tx_n)$ and, therefore, $(x_n)$ have no weakly Cauchy subsequences. Applying Rosenthal's $\ell_1$-theorem twice, we find a subsequence $(x_{n_k})$ such that $(x_{n_k})$ and $(Tx_{n_k})$ are both equivalent to the unit vector basis of~$\ell_1$. It follows that $T$ is not strictly singular.
\end{proof}

\begin{proposition}\label{SS-E-ell-p}
For all $p\in[1,\infty)$, $\mathcal{SS}(d_{w,p},\ell_p)=\mathcal K(d_{w,p},\ell_p)$.
\end{proposition}
\begin{proof}
By Proposition~\ref{SS-X-ell-1}, we only have to consider the case $p>1$. Let $T\not\in\mathcal{K}(X,\ell_p)$. Pick a bounded sequence $(x_n)$ in $X$ such that $(Tx_n)$ has no convergent subsequences. Since $d_{w,p}$ contains no copies of~$\ell_1$, by Rosenthal's $\ell_1$-theorem we may assume that $(x_n)$ is weakly Cauchy. Passing to a further subsequence, we may assume that the sequence $(Ty_n)$, where $y_n=x_{2n}-x_{2n-1}$, is seminormalized. It follows that $(y_n)$ is also seminormalized. Also, $(y_n)$ and, therefore, $(Ty_n)$ are weakly null. Passing to a subsequence of $(x_n)$, we may assume that $(y_n)$ and $(Ty_n)$ are both basic, equivalent to block sequences of $(e_n)$ and $(f_n)$, respectively. By \cite[Proposition~5]{ACL} and \cite[Proposition~2.a.1]{LT77}, $(f_n)\succeq(y_n)$ and $(f_n)\sim(Ty_n)$. So, we obtain the diagram
$$
(f_n)\succeq(y_n)\succeq(Ty_n)\sim(f_n).
$$
Hence $[y_n]$ is isomorphic to $[Ty_n]$, so that $T$ is not strictly singular. 
\end{proof}

The following lemma is standard.

\begin{lemma}\label{ell-1-dominates}
Let $X$ be a Banach space. Every seminormalized basic sequence in $X$ is dominated by the unit vector basis of~$\ell_1$.
\end{lemma}

\begin{lemma}\label{ell-1-const-perturb}
Let $(x_n)$ and $(y_n)$ be two sequences in a Banach space $X$ such that $(x_n)$ is equivalent to the unit vector basis of $\ell_1$ and $(y_n)$ is convergent. Then the sequence $(z_n)$ defined by $z_n=x_n+y_n$ has a subsequence equivalent to the unit vector basis of~$\ell_1$.
\end{lemma}
\begin{proof}
Observe that $(z_n)$ cannot have weakly Cauchy subsequences since $(x_n)$ does not have such subsequences. Since $(z_n)$ is bounded, the result follows from Rosenthal's $\ell_1$-theorem.
\end{proof}

Recall that an operator $A$ between two Banach lattices $X$ and $Y$ is called \term{positive} if $x\ge 0$ entails $Tx\ge 0$.

Conjecture~\ref{conjecture} asserts, in particular, that if $T\in\mathcal{SS}(d_{w,p})$ and $T=AB$ for some $A\colon d_{w,p}\to\ell_p$ and $B\colon\ell_p\to d_{w,p}$ then $T\in\overline{J^j}$. In the next theorem, we prove this under the additional assumptions that $p=1$ and both $A$ and $B$ are positive.

\begin{theorem}\label{positive-factor}
Let $T\in\mathcal{SS}(d_{w,1})$ be such that $T=AB$, where $A\in L(\ell_1, d_{w,1})$, $B\in L(d_{w,1}, \ell_1)$, and both $A$ and $B$ are positive. Then $T\in\overline{J^j}$.
\end{theorem}
\begin{proof}
Define a sequence $(A_N)$ of operators in $L(\ell_1, d_{w,1})$ by the following procedure. For each $n\in\mathbb N$, let $A_Nf_n$ be obtained from $Af_n$ by keeping the largest $N$ coordinates and replacing the rest of the coordinates with zeros. Since $Af_n\ge 0$ for all $n\in\mathbb N$, this defines a positive operator $\ell_1\to d_{w,1}$. Also, $\norm{A_Nf_n}\le\norm{Af_n}\le\norm A$ for all $n\in\mathbb N$, hence $\norm{A_N}\le\norm A$.

Define $A'_N=A-A_N$. It is clear that $0\le A'_Nf_n\le Af_n$ for all $n\in\mathbb N$, hence $A'_N\ge 0$ and $\norm{A'_N}\le\norm A$. We claim that $A'_{N}\to 0$ in the strong operator topology (SOT). Indeed, since $A'_Nf_n$ is obtained from $Af_n$ by removing the largest $N$ coordinates, the elements of the matrix of $A'_N$ are all smaller than $\frac{\norm{A}}{s_N}$ by Lemma~\ref{norm-entry-correspondence}. In particular, if $0\le x\in\ell_1$, then $A'_Nx\downarrow 0$; it follows that  $\norm{A'_Nx}\to 0$ because  $d_{w,1}$ has order continuous norm (see Remark~\ref{wsc-ocn}). If $x\in\ell_1$ is arbitrary then $\norm{A'_Nx}\le\bignorm{A'_N\abs{x}\,}\to 0$.

We will show that $\norm{A'_NB}\to 0$ as $N\to\infty$, so that $\norm{AB-A_NB}\to 0$ as $N\to\infty$. Since $(A_Nf_n)_{n=1}^\infty$ is almost lengthwise bounded (in fact, the vectors in the sequence $(A_Nf_n)_{n=1}^\infty$ all have at most $N$ nonzero entries), the theorem will follow from Theorem~\ref{ell-1-factor}.

Assume, by way of contradiction, that there are $c>0$ and a sequence $(N_k)$ in $\mathbb N$ such that $\norm{A'_{N_k}B}> c$. Then there exists a normalized positive sequence $(x_k)$ in $d_{w,p}$ such that $\norm{A'_{N_k}Bx_k}>c$. By Rosenthal's $\ell_1$-theorem, we may assume that $(x_k)$ is either weakly Cauchy or equivalent to~$(f_n)$.

Assume that $(x_k)$ is weakly Cauchy. Then $(Bx_k)$ is weakly Cauchy. Since $(Bx_k)$ is a sequence in~$\ell_1$, it must converge to some $z\in\ell_1$ by the Schur property. Then $\norm{A'_{N_k}Bx_k-A'_{N_k}z}\le\norm{A'_{N_k}}\cdot\norm{Bx_k-z}\le\norm{A}\cdot\norm{Bx_k-z}\to 0$. Since $A'_{N_k}\to 0$ in SOT, it follows that $A'_{N_k}Bx_k\to 0$, contrary to the assumption. Therefore $(x_k)$ must be equivalent to~$(f_n)$.

Since the entries of the matrix of $A'_N$ are all less than $\frac{\norm{A}}{s_N}$, the coordinates of the vector $A'_{N_k}Bx_k$ are all less than $\frac{\norm{A}}{s_{N_k}}\norm B\to 0$. Hence, passing to a subsequence, we may assume that $(A'_{N_k}Bx_k)$ is equivalent to a block sequence $(u_k)$ of $(e_n)$ such that each $u_k$ is a restriction of $A'_{N_k}Bx_k$. In particular, the coordinates of $(u_k)$ converge to zero. Passing to a further subsequence, we may assume by Remark~\ref{ell-p-wealth} that $(A'_{N_k}Bx_k)\sim(f_n)$.

The sequence $(Tx_{k})$ cannot have subsequences equivalent to $(f_n)$ since $T$ is strictly singular. Therefore, by Rosenthal's $\ell_1$-theorem, we may assume that $(Tx_k)$ is weakly Cauchy. Since $d_{w,1}$ is weakly sequentially complete (Remark~\ref{wsc-ocn}), the sequence $(Tx_k)$ weakly converges to a vector $y\in d_{w,1}$. Since the positive cone in a Banach lattice is weakly closed, $y\ge 0$.

Note that $Tx_k\ge A'_{N_k}Bx_{k}\ge u_k\ge 0$ for every~$k$. Since $(u_k)$ is a seminormalized block sequence of $(e_n)$, it follows that $(Tx_k)$ is not norm convergent. Write $Tx_k=y+h_k$; then $(h_k)$ converges to zero weakly but not in norm. Therefore, passing to a subsequence, we may assume that $(h_k)$ is seminormalized and basic (but not, necessarily, positive).

Let $r_k=A'_{N_k}Bx_k-(A'_{N_k}Bx_k\wedge y)\ge 0$, $k\in\mathbb N$. Observe that $A'_{N_k}Bx_k\wedge y\in[0,y]$ for all~$k$. Since $d_{w,1}$ has order continuous norm and the order in $d_{w,1}$ is defined by a 1-unconditional basis, order intervals in $d_{w,1}$ are compact (see, e.g., \cite[Theorem 6.1]{Wnuk}). Therefore, passing to a subsequence of $(x_{n_k})$, we may assume that $(A'_{N_k}Bx_k\wedge y)$ is convergent, hence, passing to a further subsequence, $(r_k)$ is equivalent to $(f_n)$ by Lemma~\ref{ell-1-const-perturb} and Theorem~\ref{small-perturbations}.

It follows from $y+h_k\ge A'_{N_k}Bx_k\ge 0$ that $\abs{h_k}\ge r_k$ for all~$k$. Passing to a subsequence, we may assume by Lemma~\ref{modify-lattice-domination} that $(h_k)\succeq(r_k)\sim(f_n)$. By Lemma~\ref{ell-1-dominates}, in fact $(h_k)\sim(f_n)$, and, hence, by Lemma~\ref{ell-1-const-perturb}, $(ABx_k)\sim(f_n)$. Since also $(x_k)\sim(f_n)$, this contradicts to $T=AB\in\mathcal{SS}(d_{w,1})$.
\end{proof}

\end{document}